\documentclass[12pt,british]{article}
\usepackage{graphicx}
\usepackage{amsmath,amsfonts,amssymb,amscd}
\usepackage{latexsym,epsfig}
\usepackage{babel}
\usepackage{theorem}
\usepackage{enumerate}
\usepackage{amsbsy}
\usepackage[matrix,arrow]{xy}
\usepackage{color}

\renewcommand{\epsilon}{\ensuremath{\varepsilon}}
\renewcommand{\phi}{\ensuremath{\varphi}}
\renewcommand{\to}{\ensuremath{\longrightarrow}}

\makeatletter

\renewcommand{\p@enumii}{}

\def\@enum@{\list{\csname label\@enumctr\endcsname}%
           {\usecounter{\@enumctr}\def\makelabel##1{
\normalfont\ignorespaces\emph{{##1}~}}
\setlength{\labelsep}{3pt}
\setlength{\parsep}{0pt}
\setlength{\itemsep}{5pt}
\setlength{\leftmargin}{0pt}
\setlength{\labelwidth}{0pt}
\setlength{\listparindent}{\parindent}
\setlength{\itemindent}{0pt}
\setlength{\topsep}{3pt plus 1pt minus 1 pt}}}
\def\@map#1#2[#3]{\mbox{$#1 \colon\thinspace #2 \to #3$}}
\def\map#1#2{\@ifnextchar [{\@map{#1}{#2}}{\@map{#1}{#2}[#2]}}
\DeclareRobustCommand*\textsubscript[1]{\@textsubscript{\selectfont#1}}
\def\@textsubscript#1{{\m@th\ensuremath{_{\mbox{\fontsize\sf@size\z@#1}}}}}
\makeatother

\DeclareMathOperator{\id}{\text{Id}}
\DeclareMathOperator{\inte}{\text{int}}

\theoremheaderfont{\normalfont\scshape}
\newtheorem{thm}{Theorem}

\newtheorem{prop}{Proposition}
\newtheorem{cor}[thm]{Corollary}

\newtheorem{que}{Question}
\theorembodyfont{\rmfamily}

\newcommand{\eop}{%
  \relax
  \ifvmode
    \noindent
  \else
    \unskip
    \hskip0pt plus-1fill\relax
  \fi
  \vrule width0pt
  \nobreak
  \hfill
  {\hspace*{\fill}\vrule width3pt height8pt depth0pt}%
}

\newenvironment{proof}{\par\vspace{\partopsep}\noindent\emph{Proof.}}
{\eop\par\vspace{\parsep}}

\newtheorem{rem}[thm]{Remark}

\def\P{\mathcal{P}}
\def\A{\mathcal{A}}
\def\X{\mathcal{E}}
\def\x{\xi}

\def\MCG{\textup{MCG}}
\def\PMCG{\textup{PMCG}}
\def\Aut{\textup{Aut}}
\def\Homeo{\textup{Homeo}}
\def\H{\textup{H}}
\def\T{\textup{T}}
\def\E{\bar{\mathcal{E}}}
\def\Hil{\textup{Hil}}
\def\P{\mathcal{P}}

\def\B{\textup{B}}
\def\M{\mathcal{M}}
\newcommand{\Sp}{\ensuremath{\mathbb S}}

\newcommand{\F}{\ensuremath{\mathbb F}}

\newcommand{\R}{\ensuremath{\mathbb R}}

\begin{document}

\title{Hilden Braid Groups}
\author{Paolo Bellingeri \and Alessia Cattabriga}

\date{\empty}
\maketitle

\begin{abstract}
Let $\H_g$ be a genus $g$ handlebody and $\MCG_{2n}(\T_g)$  be the
group of the isotopy classes of orientation preserving
homeomorphisms of $\T_g=\partial\H_g$,  fixing a given set of $2n$
points. In this paper we study two particular subgroups of $\MCG_{2n}(\T_g)$
which generalize Hilden groups defined by Hilden in \cite{Hil}. As well as Hilden groups
are
related to plate closures of braids, these generalizations are related to
Heegaard splittings of manifolds and to bridge decompositions of links. Connections between these subgroups 
and motion groups of  links in  closed 3-manifolds are also provided.

\bigskip

\end{abstract}

\section{Introduction}

In \cite{Hil} Hilden introduced and found generators for two particular subgroups of  the mapping class group of the sphere with $2n$ punctures. Roughly speaking these groups consist of  (the isotopy classes of) homeomorphisms of the puctured sphere which admit  an extension to the 3-ball fixing $n$ arcs embedded in the 3-ball and bounded by the punctures.   The interest in these groups was motivated by the theory of links in $\Sp^3$ (or in $\mathbb R^3$). In \cite{B},  Hilden's generators were used in order to find a finite number of explicit equivalence moves relating two braids having the same  plate closure. More recently, many authors investigated different groups related to Hilden's ones (see \cite{BC,BH,T,T2}) and in particular  motions groups (introduced in~\cite{G}). The second author introduced in~\cite{CM} a higher genus generalization of Hilden's groups. These groups are subgroups of punctured mapping class groups of  closed surfaces and are related  to the study of link theory in a closed 3-manifold.

In this paper we define and study a different  higher genus generalization of  Hilden's groups: the Hilden braid groups. These groups can be seen as subgroups of the ones studied in \cite{CM} and they can be thought as a generalization of  Hilden's groups in the ``braid direction''. Analogously to the genus zero case, our interest in these groups is mainly motivated by  the theory of links in a closed 3-manifold.  With respect to groups  introduced in \cite{CM}, our groups seem to be more useful on studying links in a fixed manifolds.

The paper is organized as follows: Section~2 is devoted to the definition of  Hilden braid groups and  pure Hilden braid groups;
a set of generators for these groups will be provided in  Section~\ref{second} (Theorems~2 and~3).
As we will prove in Section~\ref{first},  by fixing a Heegaard decomposition of a given 3-manifold, it is possible to   define a plat-like closure for  $2n$-string braids on the Heegaard surface. In this setting, Hilden braid groups play a role similar to   Hilden groups in the  plat closure of  classical braids.  Moreover, as in the genus zero case,  Hilden braid groups are connected with  motion groups of links in  closed 3-manifolds; this relation is established  in Section~\ref{third} (Theorem~7). 
\medskip
 
\section{Hilden groups: topological generalizations}

\begin{figure}[ht]
\label{Fig1}
\begin{center}
\includegraphics*[totalheight=6cm]{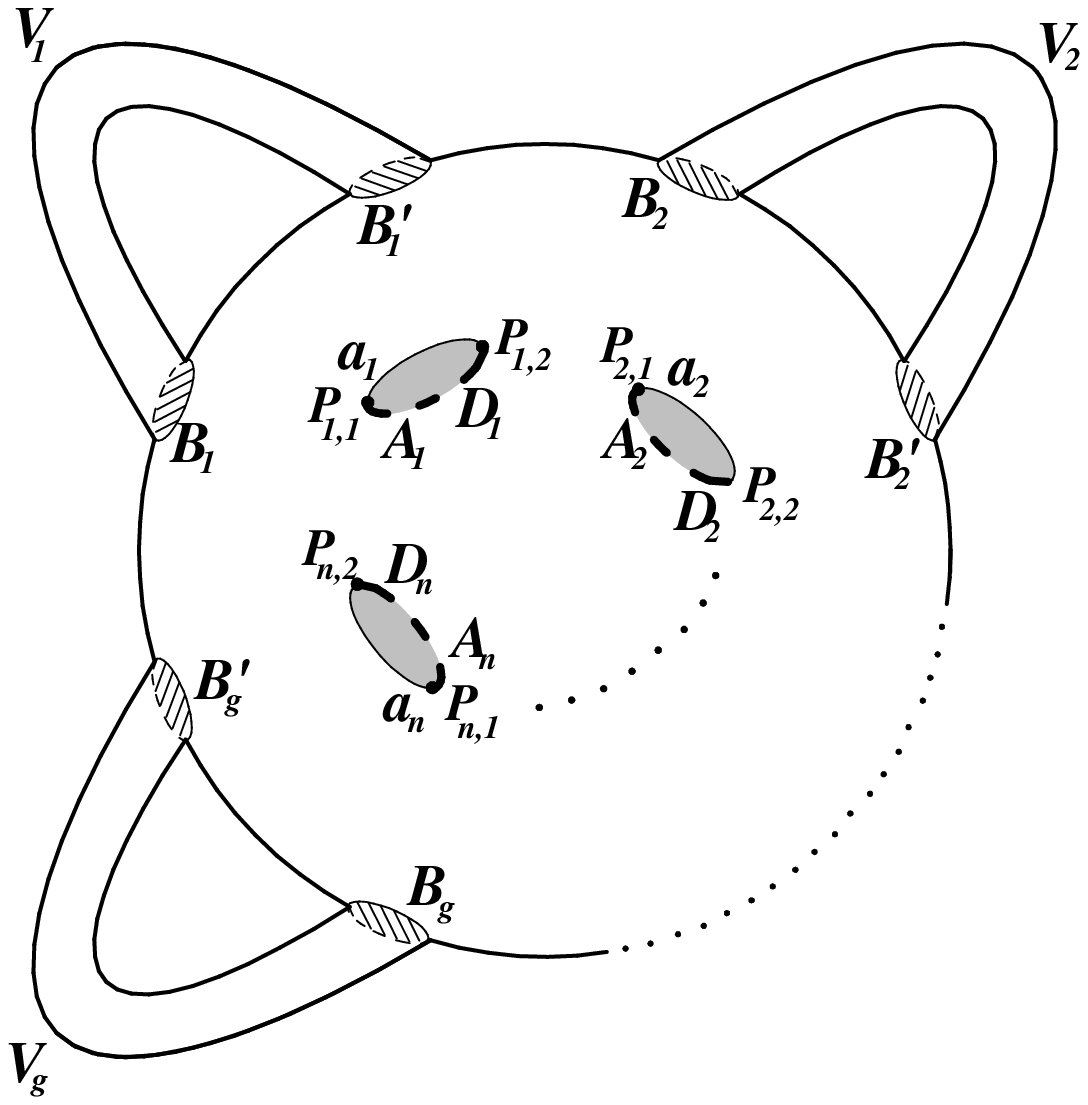}
\end{center}
\caption{The model for a genus $g$ handlebody and a trivial
system of arcs.}
\end{figure}

Referring to Figure \ref{Fig1}, let $\H_g$  be an oriented handlebody of genus $g\geq 0$ and
$\partial \H_g=\T_g$.  We recall that a
system of $n$ pairwise disjoint
arcs $\A_n=\{A_1,\ldots, A_n\}$ properly embedded in $\H_g$ is called
\textit{trivial} or \textit{boundary parallel}
if there exist $n$ disks (the grey ones in  Figure~1)
$D_1,\ldots,D_n$, called \textit{trivializing disks},
embedded in $\H_g$ such that $A_i\cap D_i=A_i\cap\partial D_i=A_i$,
$\partial D_i-A_i\subset \partial\H_g$ and $A_i\cap D_j=\emptyset$, for
$i,j=1,\ldots,n$ and $i\ne j$. 

By means of  the trivializing disks $D_i$  we can
``project'' each arc $A_i$ into the arc $a_i=\partial D_i-\inte(A_i)$ embedded in $\T_g$ and with the property  that  $a_i\cap a_j=\emptyset$, if $i\ne j$.  We
denote with $P_{i,1},P_{i,2}$  the endpoints of the arc
$A_i$ (which clearly coincide with the endpoints of  $a_i$), for
$i=1,\ldots,n$. 

Let $\MCG_{2n}(\T_g)$ (resp.
$\MCG_n(\H_g)$) be  the group of the isotopy classes of orientation
preserving homeomorphisms of $\T_g$ (resp. $\H_g$)  fixing the set
$\P_{2n}=\{P_{i,1},P_{i,2}\,\vert\, i=1,\ldots,n\}$ (resp.
$A_1\cup\cdots \cup A_n$).
  
The \textit{Hilden mapping class group} $\X_{2n}^g$  is the subgroup of
$\MCG_{2n}(\T_g)$ defined as the image of the injective group
homomorphism $\MCG_n(\H_g)\to \MCG_{2n}(\T_g)$
induced by restriction to the boundary. 
In other words, $\X_{2n}^g$
consists of the isotopy classes of homeomorphisms that admit an
 extension to $\H_g$ fixing $A_1\cup\cdots\cup A_n$.
Moreover, if $\PMCG_{2n}(\T_g)$ denotes the subgroup of
$\MCG_{2n}(\T_g)$ consisting of the isotopy classes of the
homeomorphisms of $\T_g$  fixing the punctures pointwise, we set
$\E_{2n}^g=\PMCG_{2n}(\T_g)\cap\X_{2n}^g$ and call it the \textit{
pure Hilden mapping class group}.
As recalled before, the  groups $\X_{2n}^0$ and $\E_{2n}^0$  were first introduced and studied
by Hilden in \cite{Hil}, where the author found a finite set of
generators, while in \cite{CM}  has been provided  a  set of generators for all Hilden
(pure) mapping class groups.

Now we are ready to define  Hilden braid groups. Consider the commutative diagram
$$\begin{CD}
\MCG_n(\H_g)@>\cong>>\X_{2n}^g\subset\MCG_{2n}(\T_g)\\
@VV\bar{\Omega}_{g,n}V     @VV\Omega_{g,n}V\\
\MCG(\H_g)@>\cong>>\X_0^g\subset\MCG(\T_g).
\end{CD}
$$
where the  vertical rows are forgetfull homomorphisms.
The \emph{$n$-th Hilden braid group} of the surface $\T_g$   is
the  group
$\Hil_{n}^g :=\X_{2n}^g \cap \ker \Omega_{g,n}\cong \ker \bar{\Omega}_{g,n}.$
The
\emph{$n$-th Hilden pure braid group} $P\Hil_{n}^g$ of the surface  $T_g$  is the
pure part of $\Hil_{n}^g$, that is $\Hil_{n}^g\cap\PMCG_{2n}(\T_g)$.  Notice that, since
$\MCG(\Sp^2)=\PMCG(\Sp^2)=1$ ,then $\X_{2n}^0=\Hil_{n}^0$ and
$\E_{2n}^0=P\Hil_{n}^0$. Moreover, in \cite{B2}, it is shown that $\ker(\Omega_{g,n})$ is
isomorphic to the quotient of the braid group $\B_{2n}(\T_g)$
by its center, which is  trivial if $g\geq 2$. So, if $g\geq 2$ we can see
$\Hil_{n}^g$ as a subgroup of the braid group of the surface $\T_g$.

In \cite{T,T2}, Tawn found a finite presentation for two groups that  he  called the   Hilden group $\textup{\textbf{H}}_{2n}$ and pure Hilden
group $\textup{\textbf{PH}}_{2n}$. The definition proposed by Tawn is slightly different from ours; indeed, $\textup{\textbf{H}}_{2n}$ and $\textup{\textbf{PH}}_{2n}$ are subgroups of, 
respectively,  
the braid group  $\B_{2n}$ and  the pure braid group $\textup{P}_{2n}$,
instead of, respectively, the mapping class group $\MCG_{2n}(\Sp^2)$ and the pure mapping class group $\PMCG_{2n}(\Sp^2)$. Nevertheless, from these presentations it is not difficult to obtain a presentation for $\Hil_n^0$ and $P\Hil_n^0$ as sketched in the following.
Consider the inclusion of the punctured $2n$ disk into the $2n$ punctured sphere: it 
induces  surjective maps from $\B_{2n}=\MCG_{2n}(D^2)$ to  $\MCG_{2n}(\Sp^2)$
 and   from $\textup{P}_{2n}=\PMCG_{2n}(D^2)$ to  $\PMCG_{2n}(\Sp^2)$.
The kernels of these maps coincide with the subgroup normally generated
by the center of  $\MCG_{2n}(D^2)$  and the element
$\sigma_1 \ldots \sigma_{2n-1}^2 \ldots \sigma_1$, where $\sigma_i$ denotes the usual generator of the braid group $\B_{2n}$.
One can easily show  that such  elements belong to 
$\textup{\textbf{PH}}_{2n}$ and that these maps restrict to surjective homomorphisms $\textup{\textbf{H}}_{2n}\to\Hil^0_n$ and $\textup{\textbf{PH}}_{2n}\to P\Hil^0_n$.  Therefore a finite presentation of $\Hil_{n}^0$ (respectively of
$P\Hil_{n}^0$) is given by the same set of generators of   $\textup{\textbf{H}}_{2n}$ (respectively
$\textup{\textbf{PH}}_{2n}$) and the same set of relations of   $\textup{\textbf{H}}_{2n}$ (respectively
$\textup{\textbf{PH}}_{2n}$)
plus the relation $W_1=1$ and $W_2=1$, where $W_1$ and $W_2$ are, respectively, the generator of the 
center of  $\MCG_{2n}(D^2)$  and the element
$\sigma_1 \ldots \sigma_{2n-1}^2 \ldots \sigma_1$
both written as words in the generators of $\textup{\textbf{H}}_{2n}$ (respectively $\textup{\textbf{PH}}_{2n}$). For further details on the genus zero case see \cite{CM}, while in this paper we will mainly focus  on the positive genus cases.

\section{Generators of $\Hil_{n}^g$}
\label{second}
In this section we find a set of generators for $\Hil_{n}^g$ and $P\Hil_{n}^g$.
We start by fixing some notations.

\medskip

Referring to Figure \ref{Fig1}, for each $k=1,\ldots,g$, we denote with
$V_k$  the $k$-th 1-handle (i.e.  a solid cylinder) obtained by cutting
$\H_g$ along the  two (isotopic) meridian disks $B_k$ and $B'_k$.
Moreover, we set $b_k=\partial B_k$ and $b'_k=\partial B'_k$ and call them
meridian curves. For each $i=1,\ldots,n$ the disk $D_i$ denotes the
trivializing disk for the $i$-th arc $A_i$ and $a_i=\partial D_i\setminus
\textup{int}(A_i)$.  The endpoints of both $a_i$  and $A_i$ are denoted
with $P_{i,1},P_{i,2}$ and we set $\mathcal P_{2n}=\{P_{i,1},P_{i,2}\mid
i=1,\ldots,n\}$.  We denote with $D$ a disk embedded in $\T_g$ containing
all the arcs $a_i$ and not intersecting any meridian curve $b_k$ or
$b'_k$. Finally $\delta_i$ denotes a disk in $D$ containing $a_i$ and such
that $\delta_i\cap a_j=\emptyset$ for $i\ne j$ and $j=1,\ldots,n$.

\medskip

Let us describe certain families of homeomorphisms of $\T_g$ fixing setwise $\mathcal P_{2n}$ and whose
isotopy classes belong to $\Hil_{n}^g$. We will keep the same notation for a homeomorphism and its isotopy class.
\begin{description}
\item[Intervals] For $i=1,\ldots,n$, we denote with $\iota_i$ the
homeomorphism of
$\T_g$ that exchanges the endpoints of $a_i$ inside $\delta_i$ and that is the identity outside
$\delta_i$.
The interval  $\iota_i$ is also called braid twist (see for
instance~\cite{LP}).

\item[Elementary exchanges of arcs] For $i=1,\ldots, n-1$, let $N_i$ be a
tubular neighborhood of $\delta_i\cup \beta_i\cup\delta_{i+1}$ where
$\beta_i$ is a band connecting $\delta_i$ and $\delta_{i+1}$, lying  inside $D$
and not intersecting any arc $a_j$, for $j=1,\ldots, n$. We denote with
$\lambda_i$ the homeomorphism of $\T_g$  that exchanges $a_i$ and
$a_{i+1}$, mapping $P_{i,j}$ to $P_{i+1,j}$ inside $N_i$, for $j=1,2$,  and that is the
identity outside $N_i$.

\item[Elementary twists]  We denote with $s_{i}$ the Dehn twist along the
curve $d_i=\partial\delta_i$. Notice that  $s_{i}=\iota_i^2$ in
$\MCG_{2n}(\T_g)$.

\item[Slides of arcs] Let $C$ be an oriented simple closed curve curve in
$\T_g\setminus\P_{2n}$ intersecting $a_i$ transversally in one  point.
Consider an embedded closed annulus $A(C)$ in $\T_g$  whose core is $C$,
containing $a_i$ in its interior part and such that $A(C)\cap
\P_{2n}=\{P_{i,1},P_{i,2}\}$.  We denote with $C_1$ and $C_2$ the boundary
curves of $A(C)$ with the convention that $C_1$ is the one on the left of
$C$ according to its orientation, (see Figure~3).   The slide $S_{i, C}$
of the arc $a_i$ along the curve $C$ is 
the multi-twist $T_{C_1}^{-1} T_{C_2} s_{i}^{\epsilon}$, where
$\epsilon=1$ if travelling along $C$ we see $P_{i,1}$  on the right and
$\epsilon=-1$ otherwise. Such an element fixes $a_i$ and determines on
$\T_g$ the same deformation caused by ``sliding"  the arc $a_i$ along the
curve $C$ according to its orientation. We denote the set of all the arc
slides by $\mathcal S_n^g$.

\begin{figure}[ht]
\label{arcslide}
\begin{center}
\includegraphics*[totalheight=6cm]{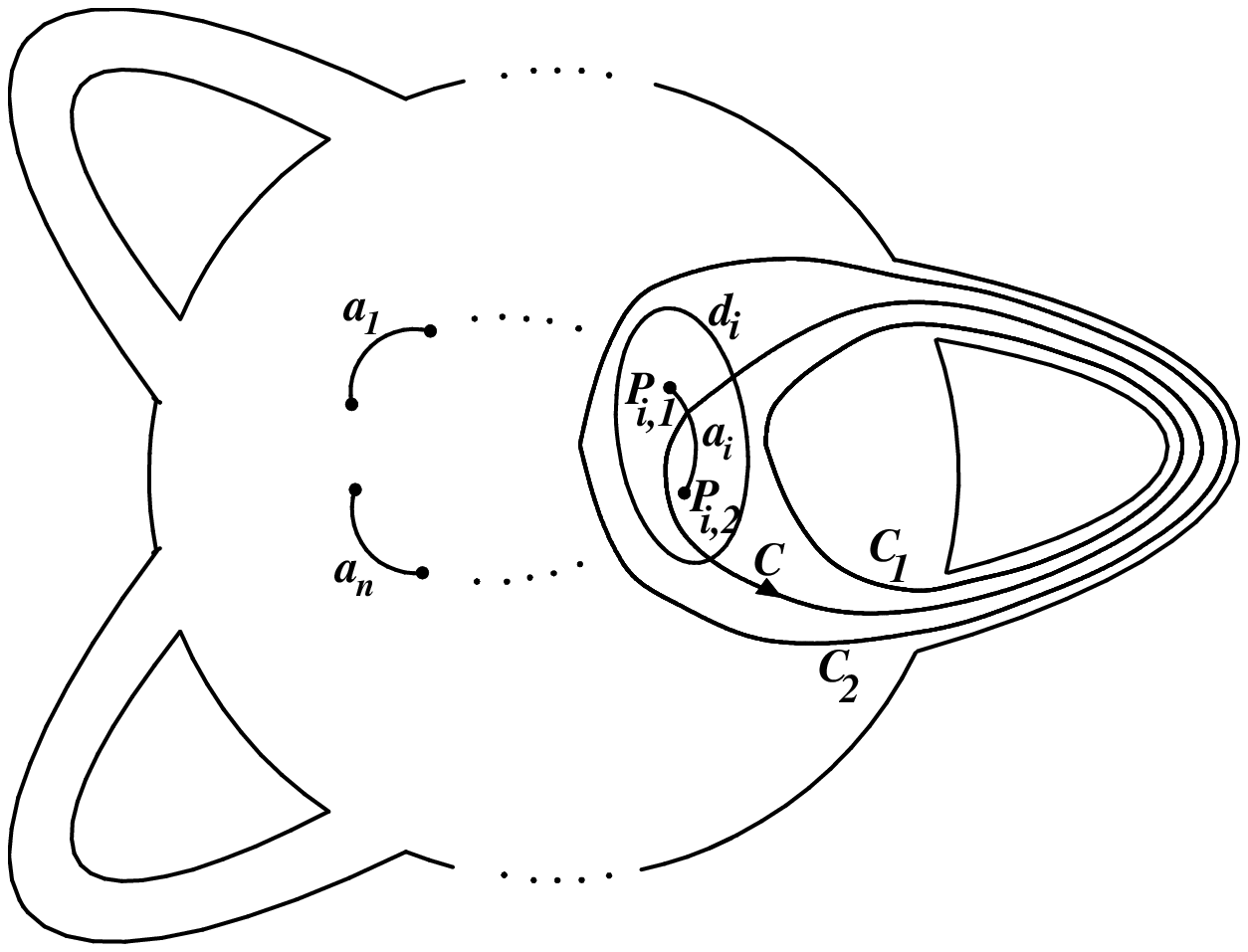}
\end{center}
\caption{The slide $S_{i,C}=T_{C_1}^{-1} T_{C_2} s_{i}$ of the arc $a_i$
along the curve $C$.}
\end{figure}

\item[Admissible slides of  meridian disks]
Let $\T_g(i)$ be the genus $g-1$ surface obtained by cutting out from
$\T_g$  the boundary of the $i$-th handle, and
capping  the resulting holes with the two meridian disks $B_i$ and $B'_i$
as in Figure~\ref{merslide}. A simple closed curve  $C$ on $\T_g(i)$ will
be called an admissible curve for the meridian disk $B_i$ if it does not
intersect $B'_i\cup\P_{2n}$, it intersect $B_i$ in a simple arc and  is
homotopic to the trivial loop in $\T_g(i)\setminus B'_i$ rel $Q$  where
$Q$ is any point of $B_i\cap C$. By exchanging the roles of $B_i$ and
$B'_i$ we obtain the definition of admissible curve for the meridian disk
$B'_i$.

Let $C$ be an admissible oriented curve for the meridian disk $B_i$.  Let
$A(C)$ be an embedded closed annulus in $\T_g(i)\setminus
(B'_i\cup\P_{2n})$  whose core is $C$ and containing $B_i$ in its
interior part.
We denote with $C_1$ and $C_2$ the boundary curves of $A(C)$ with the
convention that $C_1$ is the one on the left of $C$ according to its
orientation, (see Figure~3). 
 Notice that one between $C_1$
and $C_2$ is homotopic to $b_i$ in $\T_g(i)\setminus B'_i$, while the
other is trivial in $\T_g(i)\setminus B'_i$.   An admissible slide $M_{i, C}$ of the meridian disk $B_i$
along the curve $C$  is
the multi-twist $T_{C_1}^{-1} T_{C_2} T_{b_i}^{\epsilon}$, where
$\epsilon=1$ if  $C_1$ is homotopic to $b_i$  and $\epsilon=-1$ otherwise.
Since this homeomorphism fixes both the meridian disks   $B_i$ and $B'_i$,
it could be extended, via the identity on the boundary of $i$-th
handle, to a homeomorphism of  $\T_g$, and determines on  $\T_g$ the same
deformation caused by ``sliding"  the disk $B_i$  along the curve $C$
according to its orientation. In an analogous way we  define  $M'_{i, C}$
an admissible slide of the meridian disk $B'_i$  along an admissible
oriented curve $C$ for $B'_i$.   We denote the set of all meridian slides
with $\mathcal M_n^g$.
\end{description}

\begin{figure}[ht]
\label{merslide}
\begin{center}
\includegraphics*[totalheight=6cm]{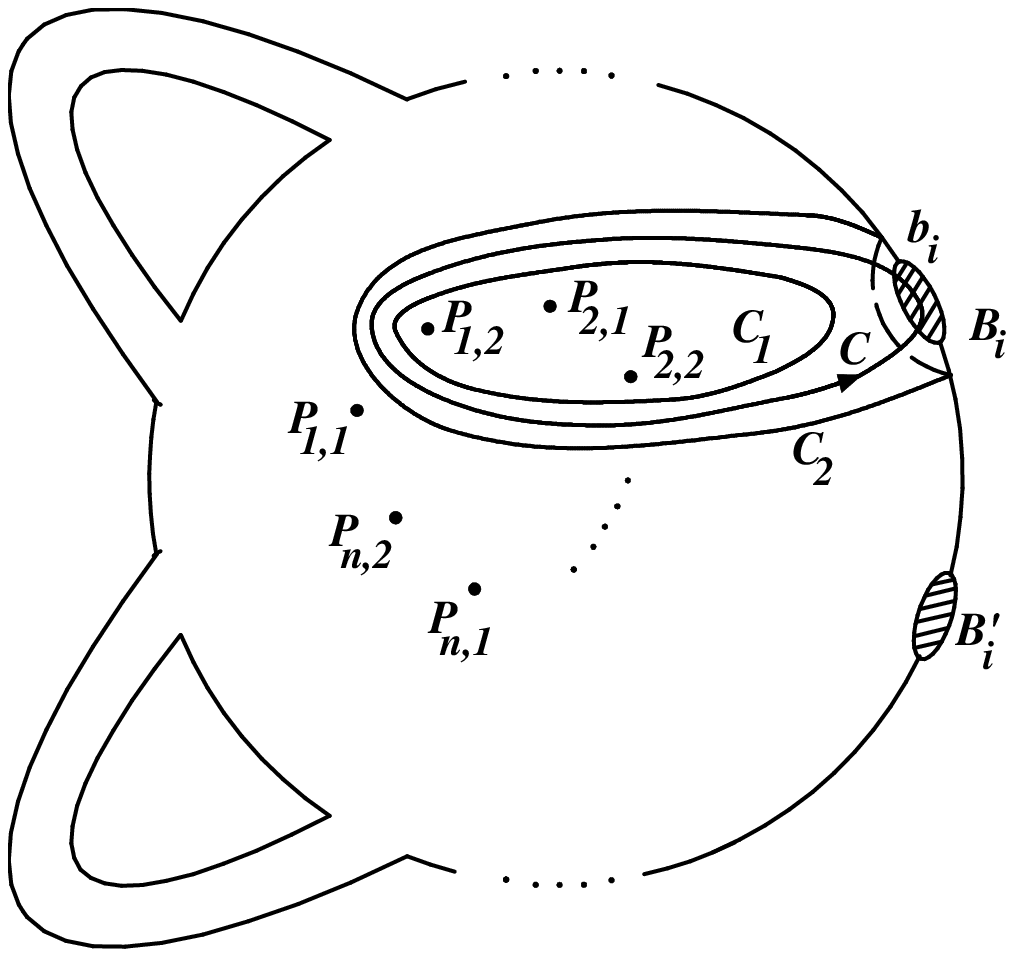}
\end{center}
\caption{The slide $M_{i,C}=T_{C_1}^{-1} T_{C_2} T_{b_i}^{-1}$ of the
meridian disk  $B_i$ along the curve $C$.}
\end{figure}

\begin{rem}
\label{slide}
In~\cite{CM}  one can find explicit extensions of all  above homeomorphisms to 
the couple $(\H_g,\mathcal A_n)$, that
is they all belongs to $\X_{2n}^g$. Moreover it is
straightforward to see that all the above elements belong also to the kernel of
$\Omega_{g,n}$ and so to $\Hil_n^g$.

It is possible to define the slide of a meridian disk $B_i$ (resp. $B_i'$)
 along a generic simple closed curve  on $\T_g(i)\setminus
B'_i\cup\P_{2n}$  (resp.  $\T_g(i)\setminus B_i\cup\P_{2n}$). Such a
meridian slide still belongs to $\X_{2n}^g$; however,  it is easy to see that a slide of a
meridian disk belong to the kernel of $\Omega_{g,n}$ (and so to $\Hil_n^g$) if and only if the
sliding curve is admissible.
\end{rem}

Let $(Z_2)^n \rtimes
 S_n$ be the signed permutation group and let
$p :\MCG_{2n}(\T_g)\to S_{2n}$ be the map which associates to any element
of $\MCG_{2n}(\T_g)$ the 
permutation induced on the punctures. 
The next proposition shows that  $\Hil_{n}^g$ is generated  by $\iota_1$,
$\lambda_i$, for $i=1,\ldots,n-1$ and  a set of generators for
$P\Hil_{n}^g$.

\begin{prop}\label{frombraidtopure}
 The exact sequence
$$1\to\PMCG_{2n}(\T_g)\to\MCG_{2n}(\T_g) \to
 S_{2n}            \to 1$$
restricts to an exact sequence $$1\to P\Hil_{n}^g \to \Hil_{n}^g \to
   (Z_2)^n \rtimes    S_n      \to 1.$$
\end{prop}
\begin{proof}
The signed permutation group can be considered as the subgroup of $S_{2n}$
generated  by the transposition
$(1\ 2)$ and the permutations $(2i-1\ 2i+1)(2i\ 2i+2)$, for
$i=1,\ldots n-1$.  Let
$\sigma\in\Hil_{n}^g$: since the extension of $\sigma$
induces a permutation of the arcs $A_1,\ldots,A_n$, then  $p(\sigma)\in
(Z_2)^n \rtimes
 S_n$. Moreover if
$p(\sigma)=1$ then an extension of it fixes the arcs pointwise, so
$\sigma\in P\Hil_{n}^g$.
\end{proof}

We say that an element $\sigma\in P\Hil_{n}^g$ is an \emph{arcs-stabilizer}
if $\sigma$ is the identity on $a_i$, for each $i=1,\ldots,n$.  The set of
all arcs-stabilizer elements of $P\Hil_n^g$  determines a subgroup of
$P\Hil_n^g$ that we call the \emph{arcs-stabilizer subgroup}.

The subgroup $FP_n(\T_g)$ of $\ker(\Omega_{g,n}) \cap \PMCG_{2n}(\T_g)$, consisting  of the elements fixing the  arcs $a_1, \ldots, a_n$
is called in~\cite{BG} \emph{$n$-th framed pure braid group of $\Sigma_g$}:
this group  is a (non trivial) generalization of pure framed braid groups
considered in~\cite{KS,MM} and several equivalent definitions have been
provided in~\cite{BG}.

\begin{prop}\label{stabilising}
Let $FP_n(\T_g)$ be the $n$-th framed braid group of $\T_g$ defined above. The
arcs-stabilizer subgroup of $P\Hil_{2n}^g$  coincides with $FP_n(\T_g)$.
In particular the arcs-stabilizer subgroup of $P\Hil_{2n}^g$  is generated by the
elementary twists $s_i$ and the slides $m_{i,j}$ and $l_{i,j}$ of the arc $a_i$ along the curves
$\mu_{i,j}$ and $\lambda_{i,j}$ depicted in 
Figure~4, for $i=1,\ldots,n$ and $j=1,\ldots,g$. 
As a consequence,  the slide $t_{i,k}$  of the arc $a_i$ along the curve $\tau_{i,k}$ depicted in Figure~5, 
is the 
composition of the above elementary twists and slides, for   $1 \le
i <k \le n$.
\end{prop}
\begin{proof}
By the above definitions, it is enough to show that $FP_n(\T_g)$ is a subgroup of $P\Hil_n^g$.
In \cite{BG} is shown that  $FP_n(\T_g)$ 
is generated by the
elementary twists $s_i$ and the  slides 
slides $m_{i,j}$ and $l_{i,j}$ of the arc $a_i$ along the curves
$\mu_{i,j}$ and $\lambda_{i,j}$ depicted in 
Figure~4, for $i=1,\ldots,n$ and $j=1,\ldots,g$. 
 Therefore $FP_n(\T_g)$ is a subgroup of $P\Hil_n^g$. 
The multitwists $t_{i,k}$, for $1 \le i <k \le n$, belong to
$FP_n(\T_g)$: in particular they
can be obtained by the above set of generators using lantern relations
(see \cite{BG}). 
\end{proof}

\begin{figure}[ht]
\label{purefr}
\begin{center}
\includegraphics*[totalheight=6cm]{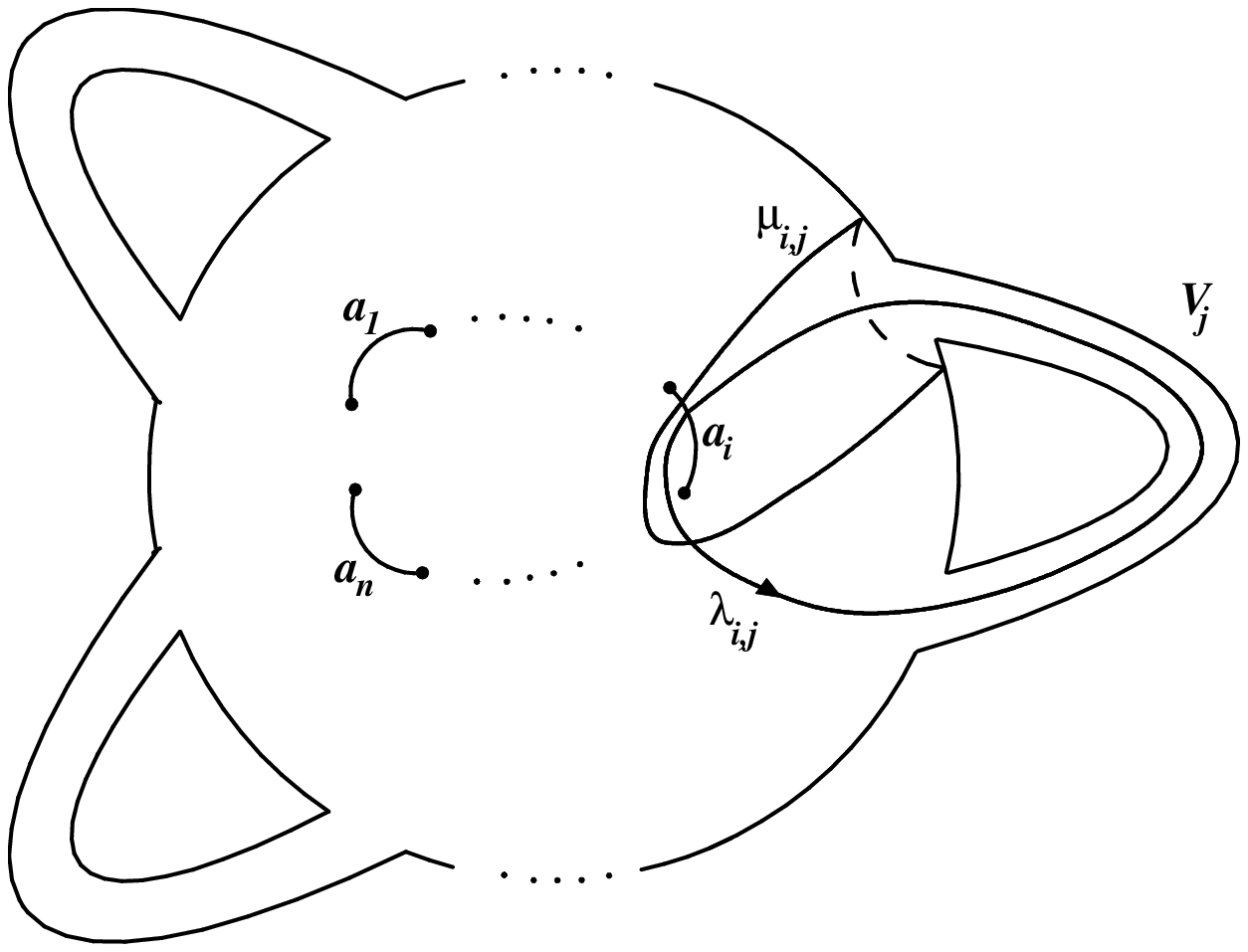}
\end{center}
\caption{The arc slides $m_{i,j}$ and $l_{i,j}$.}
\end{figure}

\begin{figure}[ht]
\label{purefr2}
\begin{center}
\includegraphics*[totalheight=6cm]{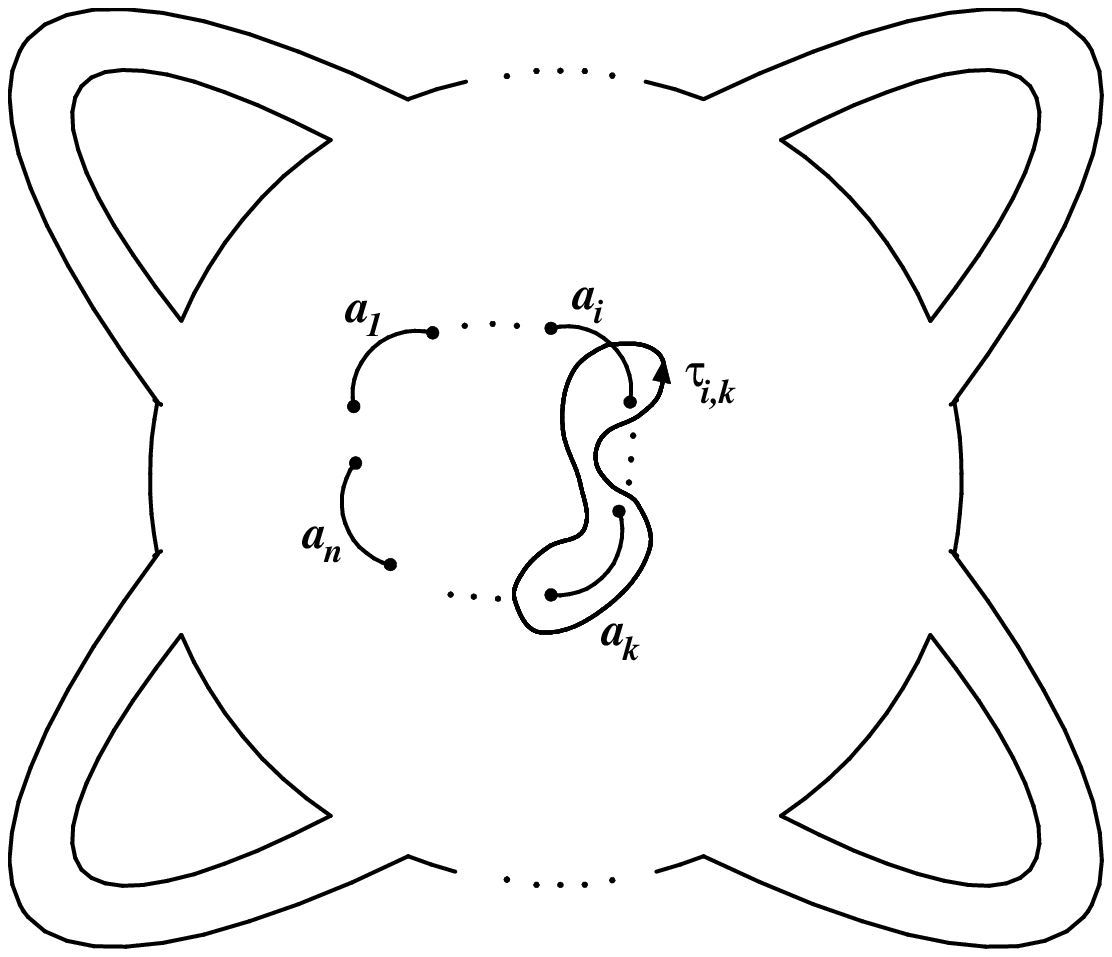}
\end{center}
\caption{The arc slide $t_{i,k}$.}
\end{figure}

Now we  will prove that an infinite set of generators for  $P\Hil_n^g$ is
given by the elementary twists, all the  arc slides and all the  admissible meridian
slides.

\begin{thm}\label{g} The group  $P\Hil_n^g$  is generated by $\mathcal
M_n^g\cup \mathcal S_n^g\cup\{s_1,\ldots,s_n\}$.
\end{thm}

\begin{proof}
Let $\mathcal{G}_n^g$ be the subgroup of  $P\Hil_n^g$  generated by
$\mathcal M_n^g\cup \mathcal S_n^g$.
By Proposition~\ref{stabilising} it is enough to show that for any
$\sigma\in P\Hil_n^g$ there exists an element $h\in\mathcal{G}_n^g$ such
that $h\sigma$ is an arcs-stabilizer. In order to do so, the first step will be to find an element $h_n\in\mathcal{G}_n^g$ such that
\begin{itemize}
\item[(i)] $h_n\sigma(a_n)=a_n$;
\item[(ii)] for each $i=1,\ldots n$ such that  $\sigma(a_n)\cap
a_i=\emptyset$ we have $h_n(a_i)=a_i$.
\end{itemize}
The element  $h_n$ will be defined as the composition of arc slides and
admissible meridian slides along opportunely chosen curves.

We denote  with $\mathcal D$ the union of all the disk $D_i$ for
$i=1,\ldots,n$ and let $\mathcal I=\sigma(D_n)\cap \mathcal D$. Up to
isotopy, we can assume that $\mathcal I$ consists of a finite number of
arcs. Clearly,  each arc $l$ in $\mathcal I$ is  a component of
$\sigma(D_n)\cap D_k$ for a unique $k$; moreover, since $\sigma(A_n)\cap
A_k=A_n\cap A_k$, if $k\ne n$ then $\sigma(A_n)\cap A_k=\emptyset$,  so
the endpoints of $l$ belong to $\sigma(a_n)\cap a_k$. If, instead, $k=n$ we
can assume that at least one endpoint of $l$ belongs to $\sigma(a_n)\cap
a_n$, since if both the endpoints lie in $A_n$ then, by composing with a
homeomorphism isotopic to the identity, the intersection $l$ can be
removed.

By an innermost argument,  it is possible to choose  $l_0\in\mathcal I$ with
$l_0\subset \sigma(D_n)\cap D_k$  such that  $l_0$
determines a disk both in $\sigma(D_n)$ and in $D_k$, whose
union is a disk $\bar D$, properly embedded in $\H_g$.
Moreover $\bar D$ intersects $\sigma(D_n)\cap \mathcal D$  only in $l_0$
if $k\ne n$, while, if  $k=n$ and one of the endpoints of $l_0$ lies in
$A_n$,  then the intersection of $\bar D$ with $\sigma(D_n)\cap \mathcal D$
  is an arc which is the union of $l_0$ with a subarc of $A_n$ going from
$l_0\cap A_n$ to one of the punctures $\partial A_n=\{P_{n,1},P_{n,2}\}$.
In any case, the boundary of $\bar D$ is the union of two simple arcs
$m_1$ and $m_2$ on $\H_g$ with $m_1\subset \sigma(a_n)$  and $m_2\subset
a_k$ (see Figure~6).

\begin{figure}[ht]
\label{proof}
\begin{center}
\includegraphics*[totalheight=6cm]{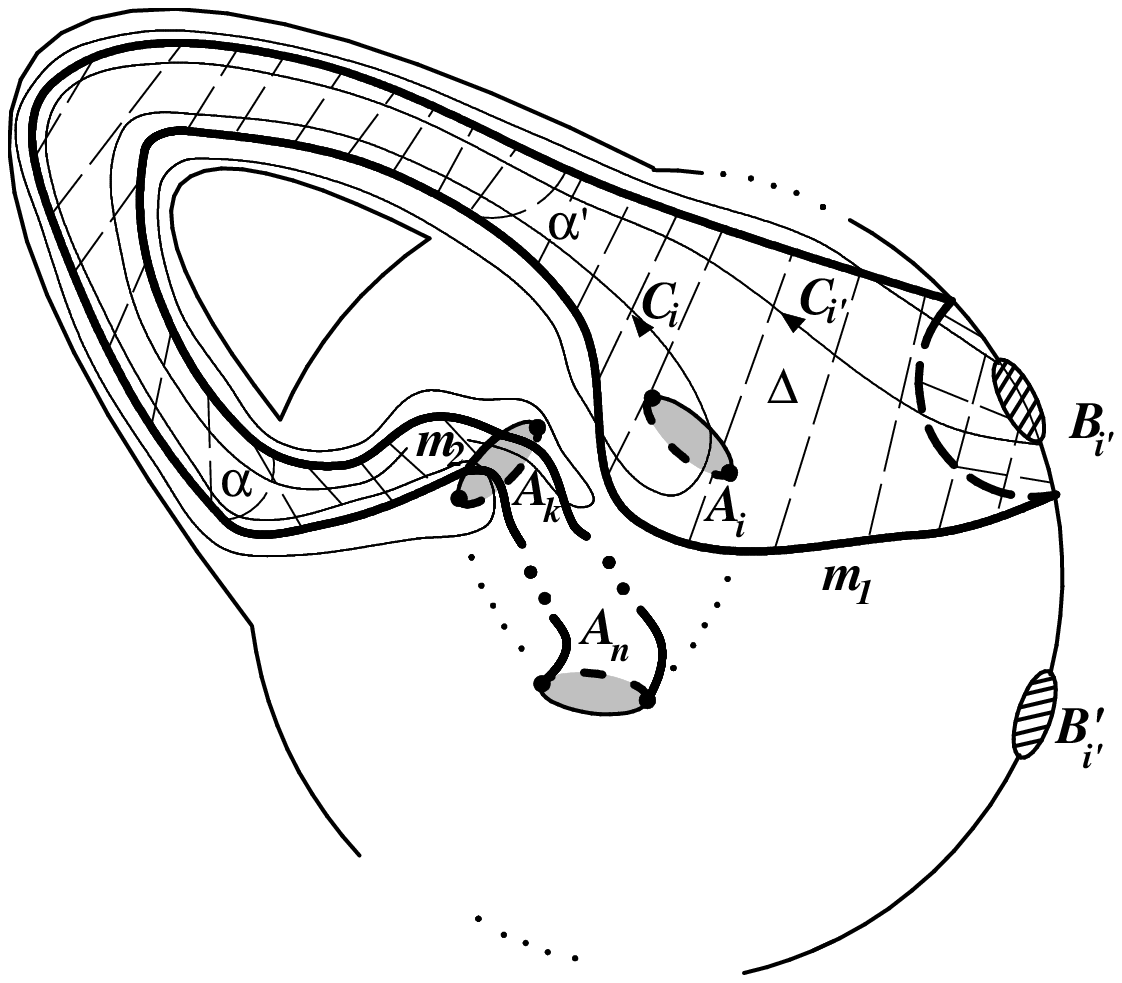}
\end{center}
\caption{}
\end{figure}

We set  $K_1=\{k\mid \bar D\cap V_k\ne\emptyset\}$ and
$K_2=\{1,\ldots g\}\setminus K_1$. Notice that if $k\in K_1$, then for
each arc $\alpha\in\bar D\cap B_k$, there exists a corresponding  arc
$\alpha'\in\bar D\cap B'_k$ such that  the union of the disks bounded by
$\alpha$ and $\alpha'$ on $B_k$, $B_k'$ and $\bar D$ is a properly
embedded disk in the handle $V_k$ bounding a ball (see Figure~6).

Indeed, if this is not the case,  the intersection $\alpha$
could be removed composing with an element isotopic to the identity.  If
$\H_g(K_2)$  denotes the handlebody (of genus $g-\vert(K_2)\vert$)
obtained from $\H_g$ by removing all the handles  $V_k$ with $k\in K_2$,
then  $\bar D$ separates   $\H_g(K_2)$   into two  connected components
$\Delta_1$ and $\Delta_2$. Let $\Delta$
be the the connected component of  $\H_g(K_2)\setminus \bar D$ that
does not contain $A_k$, which is the dotted zone in Figure \ref{proof}.

Each disk $D_i$ with $i\ne k$ and each meridian disk $B_k$ or $B_k'$, with
$k\in K_2$ is either contained or disjoint from $\Delta$. Let $I_1=\{i\mid
A_i\subset \Delta\}$, $I_2=\{i\mid i\in K_2, B_i\subset \Delta\}$ and
$I_3=\{i\mid i\in K_2, B'_i\subset \Delta\}$.  For each $i\in I_1\cup
I_2\cup I_3$  we choose a simple closed oriented curve  $C_i$ on
$\H_g(K_2)$    such that

\begin{quote}
(*)  both $C_i\cap h(a_n)=C_i\cap m_1$ and   $C_i\cap a_k=C_i\cap m_2$
consist of a single point; travelling along $C_i$ the intersection with
$a_k$ comes before the one with $\sigma(a_n)$; if $i\in I_1$ then $C_i$
intersects $a_i$ in a single point,  while if $i\in I_2$ (resp. $I_3$) then
$C_i$ is an admissible curve for the meridian disk  $B_i$ (resp. $B'_i$);
$C_i$ does not intersects all the others arcs $a_j$ and all the other
meridian disks $B_j,B'_j$, with $j\in K_2$.
\end{quote}
To see that such a $C_i$ exists, consider an arc $a$ that starts from
$a_i$, $B_i$ or $B'_i$ travels inside $\partial \Delta$ to $m_1$ without
intersecting all the other arcs $a_j$ and all the other meridian disks
$B_j,B'_j$, goes along $m_1$ to $a_k$ and finally goes along $a_k$ to an
endpoint of $a_k$ (that is $P_{k,1}$ or $P_{k,2}$) without travelling
along $m_2$. Then $C_i$ can be chosen as the boundary of a small tubular
neighboorhood of $a$. In  Figure \ref{proof} it is depicted the case
$I_1=\{i\}$, $I_1=\{i'\}$ and $I_3=\emptyset$.

If we set $f_0=\prod_{i\in I_1} S_{i, C_i} \prod_{i\in
I_2}M_{i,C_i}\prod_{i\in I_3}M'_{i,C_i}$ then $f_0(\Delta)$ does not
contain any arc and any meridian disk. This means that $f_0(\Delta)$
bounds a ball in $\H_g\setminus \A_n$ and so, up to composing with an
element isotopic to the identity,  we can remove the intersection $l_0$.
Moreover, $f_0\sigma(D_n)\cap  \mathcal D=\mathcal I\setminus l_0$.

By repeating the above procedure a finite number of times, we get an
element $f=f_kf_{k-1}\cdots f_0\in\mathcal{G}_n^g$  such that
$f\sigma(D_n)\cap \mathcal D=A_n$.  Then $f\sigma(D_n)\cup D_n$ is a
properly embedded disk $\widetilde D$ in $\H_g$. If we denote with
$\H_g(K_2)$   the handlebody obtained  from $\H_g$ by removing all the
handles that have no intersection with $\tilde D$ (that is $K_2=\{k\mid
V_k\cap \widetilde D=\emptyset\}$), then $\widetilde D$ separates
$\H_g(K_2)$   into two connected components  $\Delta_1$ and $\Delta_2$. As
above, each disk $D_i$ with $i\ne n$ and each meridian disk $B_i$ or
$B_i'$, corresponding to the removed handles, is either contained or
disjoint from $\Delta_k$, for $k=1,2$. If there exists $k=1,2$ such that
$\Delta_k$ does not contain any  disk $D_i,B_i$ or $B'_i$ then, up to composing
with an element isotopic to the identity,  we have $f\sigma(D_n)=D_n$ and
so $h_n=f$. If this is not the case, we choose one of the two connected
components, for example $\Delta_1$,  and, as before, we set
$I_1=\{i\mid i\in K_2, A_i\subset \Delta_1\}$, $I_2=\{i\mid i\in K_2, B_i\subset
\Delta_1\}$ and $I_3=\{i\mid i\in K_2, B'_i\subset \Delta_1\}$ and for
each $i\in I_1\cup I_2\cup I_3$  we choose a simple oriented closed  curve
 $C_i$ on $\H_g(K_2)$ satisfying (*). Then by taking  $h_n=\prod_{i\in
I_1} S_{i, C_i} \prod_{i\in I_2}M_{i,C_i}\prod_{i\in I_3}M'_{i,C_i}f$ we
get $h_n\sigma(a_n)=a_n$. Moreover $h_n$ satisfies (ii) since it is the compositions of slides (of
arcs or of meridian disks) along curves that by (*) intersects only the
arcs $a_i$ such that $a_i\cap \sigma(a_n)\ne\emptyset$, and so fixes all
the other arcs.

We can repeat the same procedure on $a_{n-1}$, that is we can find
$h_{n-1}\in\mathcal{G}_n^g$ with $h_{n-1}h_n\sigma(a_{n-1})=a_{n-1}$ and
such that for each $i=1,\ldots, n$ such that  $h_n\sigma(a_{n-1})\cap
a_i=\emptyset$ we have $h_{n-1}(a_i)=a_i$. This implies that
$h_{n-1}h_n\sigma(a_n)=a_n$, that is $h_{n-1}h_n\sigma$ fixes the last two
arcs. Proceeding in this way, we construct, for each $i=1,\ldots,n$, an element
$h_i\in\mathcal{G}_n^g$ such that $h_ih_{i+1}\cdots h_n\sigma(a_j)=a_j$
for each $j\geq i$. So  $h=h_1h_2\cdots h_{n-1}h_n$ is the required
element of  $\mathcal{G}_n^g$, that is $h\sigma$ is arc-stabilizer.
\end{proof}

In order to find a finite set of generators it would be enough to show that 
the subgroup of $P\Hil_n^g$ generated by $\mathcal S_n^g$ 
and the one  generated by  the $\mathcal M_n^g$ are finitely generated. The next two propositions show that the first
subgroup is finitely generated, and  the second one is  finitely
generated when $g=1$.

\begin{prop}\label{proparcslide}The subgroup of  $P\Hil_n(\T_g)$ generated by $\mathcal
S_n^g$  is finitely generated by
\begin{itemize}
\item[(1)] the elementary twist $s_i$, for ,$i=1,\ldots,n$;  
\item[(2)] the slides $m_{i,j}$ and $l_{i,j}$ of the arc $a_i$ along the curves
$\mu_{i,j}$ and $\lambda_{i,j}$ depicted in 
Figure~4, for $i=1,\ldots,n$ and $j=1,\ldots,g$.  
\item[(3)]  the slides  $s_{i,k}$, of the arc $a_i$ along the curve $\sigma_{i,k}$ depicted in Figure~7,  for $1\le i\not=k \le n$.
\end{itemize}
\end{prop}

\begin{proof}
By the definition of slide if  $C\simeq C_1\cdots C_r$ in
$\pi_1(\T_g\setminus\{P_{j,k}\mid j=1,\ldots ,n, j\ne i, k=1,2\}, *)$ then
$S_{i,C}=S_{i,C_r}\cdots S_{i,C_1}$ in $\MCG_{2n}(\T_g)$, where $*\in a_i$
is any fixed point. Then all the slides of the $i$-th arc are generated by
the slides of the $i$-th arc along a set of generators for
$\pi_1(\T_g\setminus\{P_{j,k}\mid j=1,\ldots ,n, j\ne i, k=1,2\}, *)$. The
set of slides $m_{i,j}$, $l_{i,j}$
for   $j=1,\ldots,g$ and $s_{i,k}, s'_{i,k}$, for $k=1,\ldots,n$, $i\not=k$ is therefore a
set of generators for all the slides of the $i$-th arc, where $s'_{i,k}$ is the slide of the arc $a_i$ along the curve $\sigma'_{ik}$ depicted in Figure~7.
Applying a lantern relation one obtains that $s_{i,k} s'_{i,k}=t_{i,k}
s_i^{-1}$. The statement therefore, 
follows from the fact that, by Proposition~\ref{stabilising}, the slide $t_{i,k}$ can be written as composition
of  the slides  $m_{j,r}$, $l_{j,r}$ and elementary twists.
\end{proof}

\begin{figure}[ht]
\label{slidesij}
\begin{center}
\includegraphics*[totalheight=6cm]{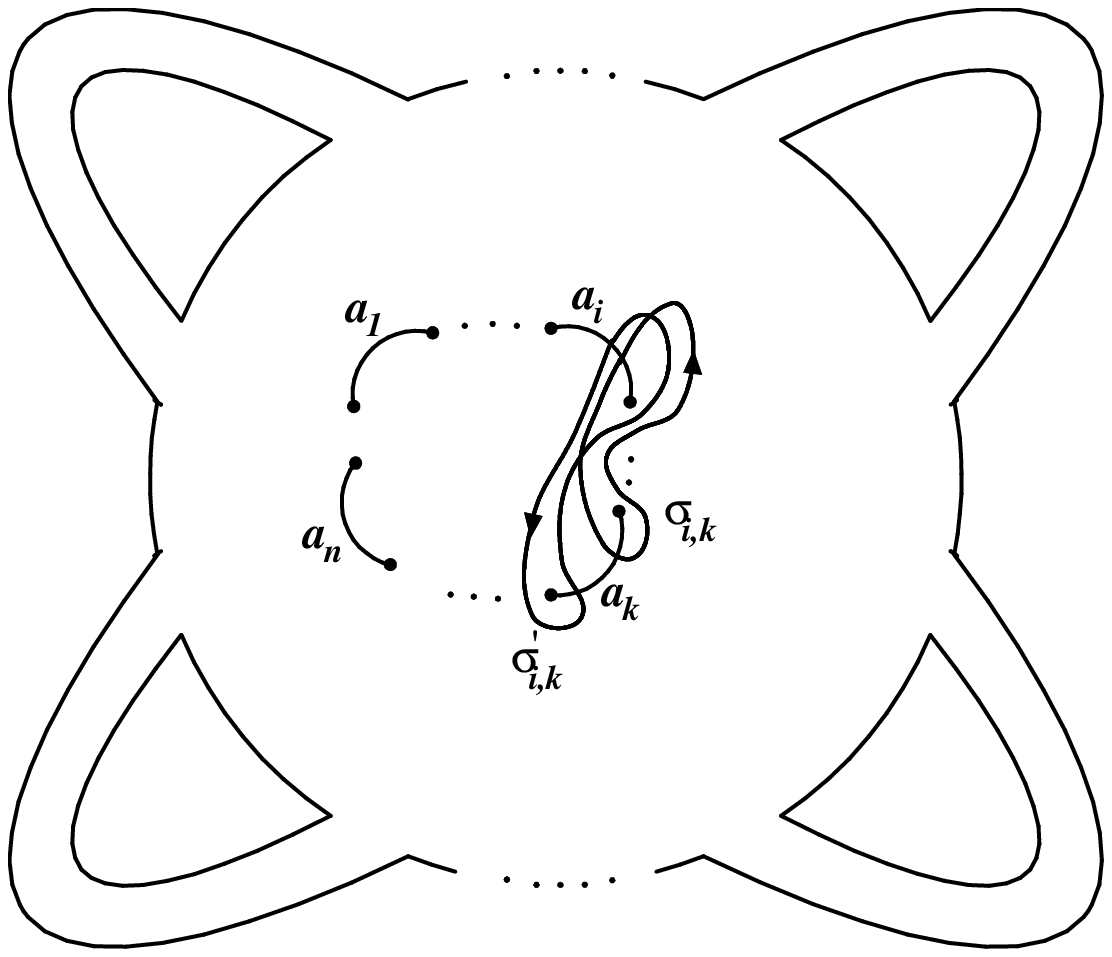}
\end{center}
\caption{The arc slides $s_{i,k}$ and $s'_{i,k}$.}
\end{figure}

\begin{prop}
The subgroup of  $P\Hil_n(\T_1)$ generated by  $\mathcal M_n^1$    is
finitely generated. 
\end{prop}

\begin{proof}
Let $\Sp^2$ be the sphere obtained by cutting out from $\T_1$  the boundary
of the handle, and
capping  the resulting holes with the two meridian disks $B_1$ and $B'_1$
as in the definiton of meridian slides. Any simple closed curve  $C$ on
$\Sp^2$
which does not intersect $B'_1\cup\P_{2n}$ and  intersects $B_1$ in a
simple arc, is an admissible curve  for the meridian disk $B_1$. 
Therefore,  if  we write $C$ as a product $C_1\cdots C_r$ of
a finite set of generators
of $\pi_1(\Sp^2\setminus (\mathcal P\cup B'_1), *)$, where $*\in B_1$ is any
fixed point, we have $M_{i,C}=M_{i,C_r}\cdots M_{i,C_1}$.  An analogous statement  holds
for $M'_{i,C}$.
\end{proof}

From these Propositions it follows that $P\Hil_n^1$ is finitely generated.
The problem of whether $P\Hil_n^g$ is finitely generated or not while $g>1$
remains open.
We end this section by giving the generators for $\Hil_n^g$.

\begin{thm}
The group $\Hil_n^g$ is generated by  
\begin{itemize}
\item[(1)] $\iota_1$ and   $\lambda_j$,  $j=1,\ldots n-1$;
\item[(2)] $m_{1,k}$, $l_{1,k}$, $s_{1,r}$ with
$k=1,\ldots,g$ and  $r=1,\ldots n$;
\item[(3)] the elements of $\mathcal M_n^g$.
\end{itemize}
\end{thm}

\begin{proof}
 The statement follows from Theorem \ref{g}, Proposition \ref{proparcslide} and  the remark that
 $s_1=\iota_1^2$ and that arc slides of the arc $a_i$ can be  reduced to arc slides of the arc $a_1$
by using (compositions of) elementary exchanges of arcs.
\end{proof}


\section{Generalized plat closure}\label{first}
One of the main motivations  to study topological generalizations of Hilden
groups comes from
link theory in 3-manifolds. In this section  we  describe a representation
of links in 3-manifolds via braids on closed
surfaces: this approach generalizes the concept of plat closure
and  explains the
role played by  $\Hil_n^g$ in this representation. We start by recalling the definition of $(g,n)$-links.

\medskip

 Let $L$ be a link in a 3-manifold $M$. We say that $L$ is a $(g,n)$\textit{-link} if
there exists a genus $g$  Heegaard surface $S$ for $M$ such that
\begin{itemize}
\item[(i)] $L$ intersects $S$ transversally and
\item[(ii)]the intersection of $L$ with both of  handlebodies
into which $M$ is divided by $S$, is a trivial system of $n$   arcs.
\end{itemize}
Such a decomposition for $L$ is called $(g,n)$\textit{-decomposition} or $n$-\textit{bridge
decomposition of genus} $g$.    The minimum $n$ such that
$L$ admits a $(g,n)$-decomposition is called \textit{genus } $g$ \textit{bridge number} of $L$.

Clearly if $g=0$ we get the usual notion of bridge decomposition and
bridge number
of links in the 3-sphere (or in $\R^3$).  Given two links $L\subset M$ and
$L'\subset M'$
we say that $L$ and $L'$ are \textit{equivalent}
if there exists an orientation preserving homeomorphism $f:M\to M'$ such
that $f(L)=L'$ and we write $L\cong L'$.

The notion of $(g,n)$-decompositions was used in \cite{CM} to develop an
algebraic representation of $\mathcal L_{g,n}$, the set of equivalence
classes  of   $(g,b)$-links, as follows. Let $(\H_g,\A_n)$ be as in Figure~1
and let $(\bar\H_g,\bar\A_n)$ be a homeomorphic copy of
$(\H_g,\A_n)$. Fix an orientation reversing homeomorphism
$\tau:\H_g\to\bar \H_g$ such that $\tau(A_i)=\bar A_i$, for each
$i=1,\ldots,n$. Then the following  application is well defined and
surjective
\begin{equation}
\label{eqrepr}
\Theta_{g,n}:\MCG_{2n}(\T_g)\longrightarrow\mathcal L_{g,n}\ \ \ \
\Theta_{g,n}(\psi)=L_{\psi}
\end{equation}
where  $L_{\psi}$ is the $(g,n)$-link in the 3-manifold $M_{\psi}$
defined by $$(M_{\psi},L_{\psi})=(\H_g,\mathcal
A_n)\cup_{\tau\psi}(\bar{\H}_g, \bar{\A}_n).$$
This means that it is possible to describe each link admitting a
$(g,n)$-decomposition in a certain 3-manifold by a
 element of $\MCG_{2n}(\T_g)$.  This element is not unique, since we have
the following result.
\begin{prop}[\cite{CM}]\label{prop1}
 If $\psi$ and $\psi'$ belong to the same double coset of $\X_{2n}^g$ in
$\MCG_{2n}(\T_g)$ then $L_{\psi}\cong L_{\psi'}$.
\end{prop}
Therefore,  in order to describe all $(g,n)$-links via \eqref{eqrepr} it
is enough to consider   double coset classes of $\X_{2n}^g$ in
$\MCG_{2n}(\T_g)$.
   This  representation has revealed to be a useful tool for studying
links in 3-manifolds, see \cite{C,CM2,CMV,K}.
However, if we represent links using  \eqref{eqrepr},  we have to deal
with links that lie in different manifolds.  If we want to  fix the
ambient manifold,
then the following remark holds.

\begin{rem}
\label{l1}
If $\psi_1,\psi_2\in\MCG_{2n}(\T_g)$ are such that
$\Omega_{g,n}(\psi_1)=\Omega_{g,n}(\psi_2)$ then $L_{\psi_1}$ and
$L_{\psi_2}$ belong to the same ambient manifold.
\end{rem}
So, in order to fix the ambient manifold,  we want to modify
representation \eqref{eqrepr} by separating  the part that determines the
manifold from the part that determines the link.

Referring to Figure~1, let $D$ be  a disk embedded in $\T_g$ containing
all the arcs $a_i$ and not intersecting any meridian curve $b_k$ or
$b'_k$, for $i=1,\ldots,n$ and $k=1,\ldots,g$.  Let $\mathcal{T}_n^g$ be
the subgroup of $\MCG_{2n}(\T_g)$ generated by Dehn twist along curves
that do not intersect the disk $D$.  We have the following proposition.
\begin{prop}
For each $\psi\in\mathcal T_{n}^g$  the link $L_{\psi}$ is a $n$ component
trivial link in $M_{\psi}$.
\end{prop}
\begin{proof}
Since, the action of $\psi$ on the punctures is trivial, for  each
$i=1,\ldots ,n$, the  arc $A_i$ is glued, via $\tau\psi$, to $\bar A_i$,
giving rise to  a connected component  of $L_{\psi}$. Moreover for each
$i=1,\ldots,n$ we have $\psi(a_i)=a_i$, so if we  set $\tau(D_i)=\bar
D_i\subset\bar\H_g$, the $i$-th  component $A_i\cup_{\tau\psi}\bar A_i$ of
 $L_{\psi}$ bounds in $M_{\psi}$  the embedded disk
$D_i\cup_{\tau\psi}\bar D_i$.
\end{proof}

Let $\T_{g,1}$ be the compact surface obtained by removing the interior
part of the disk $D$. The natural inclusion of
$\T_{g,1}$  into $\T_{g}$ with $2n$ marked points induces an injective map
$\MCG(\T_{g,1})\to \MCG_{2n}(\T_g)$ (see~\cite{RP}) and the mapping class
group
 $\MCG(\T_{g,1})$ turns out to be isomorphic to the group $\mathcal T_{n}^g$.
 On the other hand, we have also the following exact sequence:
$$1\to \pi_1(U\T_{g,1})\to\mathcal T_n^g\to\MCG(\T_g)\to 1$$
where $U\T_{g,1}$ is  the unit tangent bundle of $\T_{g,1}$ (see
\cite{BG}). As a consequence each element of $\MCG(\T_g)$ admits a 
lifting as an element of $\mathcal T_n^g$, so we  can realize any genus 
$g$ Heegaard decomposition of a
3-manifold  $M$ using an  element of
 $\mathcal T_n^g$, for any $n>0$.
Now we are ready to define the \emph{generalized plat closure}. Let  $M$
be a fixed manifold,  and  choose an element $\psi\in\mathcal T_n^{g}$
such that  $M=M_{\psi}$.   We define a map
\begin{equation}
\label{repr2}
\Theta_{g,n}^{\psi}:\ker(\Omega_{g,n})
\longrightarrow\{(g,n)-\textup{links in } M_{\psi} \}
\end{equation}
given by $\Theta_{g,n}^{\psi}(\sigma)=\Theta_{g,n}(\psi\sigma)$. We set
$\hat{\sigma}^{\psi}=\Theta_{g,n}^{\psi}(\sigma)$.

\begin{rem}
 As recalled  before, $\ker(\Omega_{g,n})$ is isomorphic to the braid
group $\B_{2n}(\T_g)$, quotiented by its center, which is trivial if
$g\geq 2$.  This means that we can interpretate $\Theta_{g,n}^{\psi}$ as
a generalization of the notion of plat closure for classical braids, as
shown schematically  in  Figure~8.
 Indeed, for  $g=0$ and $\psi=\textup{id}$ we obtain the classical plat closure.
This is the only representation in the case of the 3-sphere  (i. e. with
$g=0$),  since  $\mathcal T_n^0$ is trivial. On the contrary, the
generalized plat closure in a 3-manifold different from $\Sp^3$ 
depends  on the choice of the element $\psi\in\mathcal T_n^g$, and,
topologically, this  corresponds to the choice of a Heegaard surface of
genus $g$  for $M$.
\end{rem}

\begin{figure}[ht]
\label{plat}
\begin{center}
\includegraphics*[totalheight=6cm]{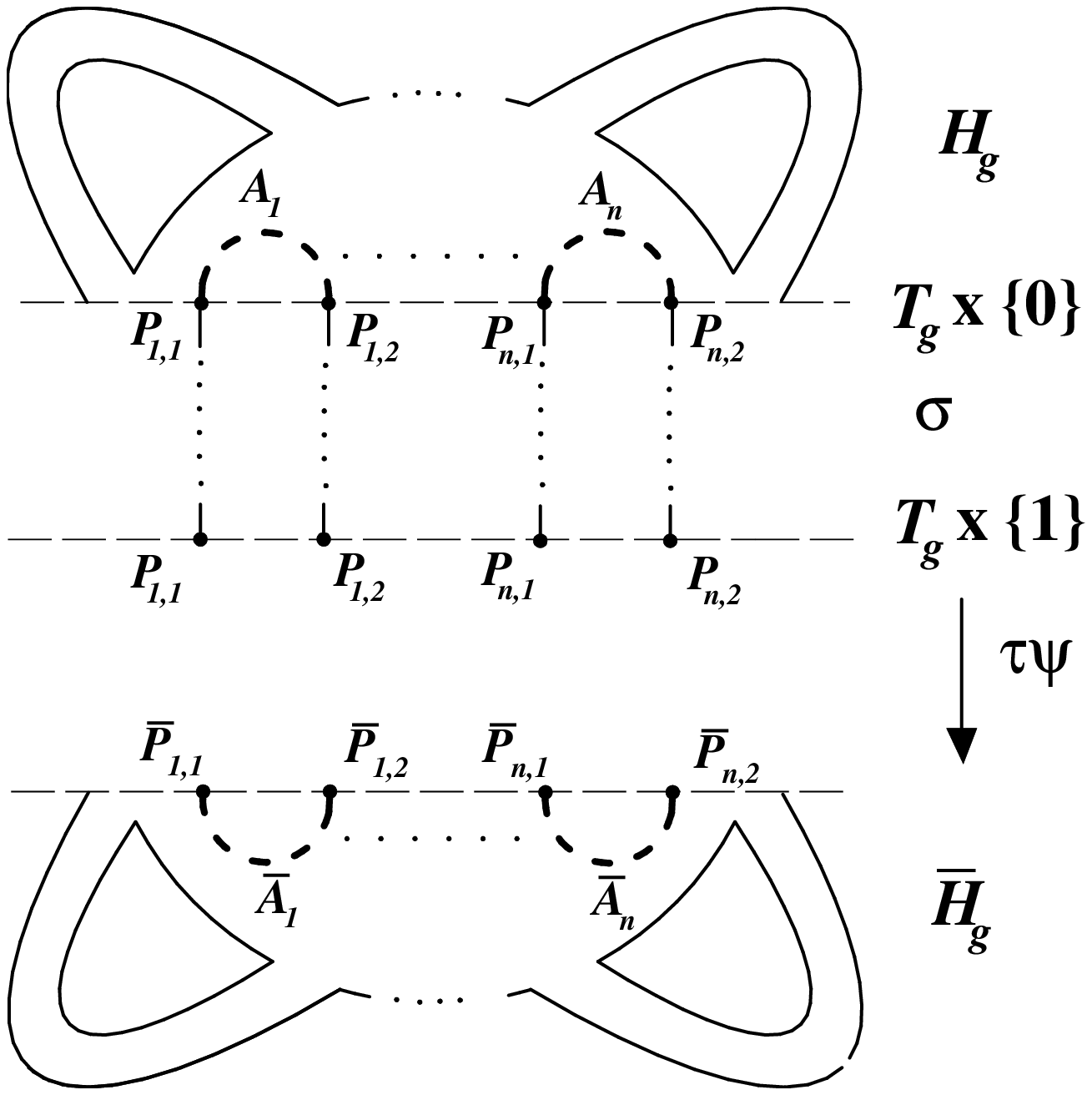}
\end{center}
\caption{A generalized plat closure}
\end{figure}

In this setting a natural question arises.
\begin{que}
\label{p1}
Is it possible to determine  when two element
$\sigma_1\in\ker(\Omega_{g,n_1})$ and $\sigma_2\in\ker(\Omega_{g,n_2})$
determine equivalent links via \eqref{repr2}?
\end{que}
A partial answer  is given by the following
statement, that is a straightforward corollary of Proposition~\ref{prop1}.

\begin{cor}
Let  $\psi\in\mathcal T_n^{g}$. Denote with
$\Hil_n^g(\psi)=\psi^{-1}\Hil_n^g\psi$.
\begin{itemize}
 \item[1)]  if $\sigma_1$ and $\sigma_2$ belong to the same  left  coset
of $\Hil_{n}^g$ in $\ker(\Omega_{g,n})$ then  $\hat{\sigma_1}^{\psi}$ and
 $\hat{\sigma_2}^{\psi}$ are equivalent links in the manifold $M_{\psi}$.
\item[2)]  if $\sigma_1$ and $\sigma_2$ belong to the same  right  coset
of $\Hil_{n}^g(\psi)$ in $\ker(\Omega_{g,n})$ then $\hat{\sigma_1}^{\psi}$
and  $\hat{\sigma_2}^{\psi}$  are equivalent links in the manifold
$M_{\psi}$.
\end{itemize}
\end{cor}
\begin{proof}
To prove the equivalence it is enough to exhibit two orientation
preserving homeomorphisms $f:(\H_g,\A_n)\to (\H_g,\A_n)$ and $\bar f:
(\bar{H}_g, \bar{\A}_n)\to (\bar{H}_g, \bar{\A}_n) $  making the following
diagram commute
$$\begin{CD}
(\partial \H_g,\partial \A_n)@>\tau\psi\sigma_1>>(\partial
\bar{\H}_g,\partial \bar{\A}_n)\\
@VVf_{\vert_{\partial}}V     @V\bar{f}_{\vert_{\partial}}VV\\
(\partial \H_g,\partial \A_n)@>\tau\psi\sigma_2>>(\partial
\bar{\H}_g,\partial \bar{\A}_n)
\end{CD}
$$
In the first case there exists $\epsilon\in\Hil_n^g$ such that
$\sigma_2=\sigma_1\epsilon$ so we can choose $f=\epsilon^{-1}$ and $\bar
f=\textup{id}$.
In the second case there exists  $\epsilon\in\Hil_n^g$ such that
$\sigma_2=\psi^{-1}\epsilon\psi\sigma_1$ so we can choose $f=\textup{id}$
and $\bar f=\tau\epsilon\tau^{-1}$.
\end{proof}

In the case of classical plat closure Question~\ref{p1} was solved  in
\cite{B}, where it is shown that two braids determine  the same plat
closure if and only if they are related by a finite sequence of moves
corresponding to generators of  $\Hil^0_n$ and a  stabilization move (see also
\cite{Mo}). 
\medskip

Another non trivial question concerns the surjectivity of the map
\eqref{repr2}. The following Proposition deals with  this problem.
\begin{prop}
\label{standard} Let $M$ be a 3-manifold with a finite number of
equivalence classes of genus $g$ Heegaard
splittings\footnote{Two Heegaard splittings of a manifold $M$
are said equivalent if there exists an homeomorphism $f:M\to M$
that send the one splitting surface into the other.}.
Then there exist  $\psi_1,\cdots, \psi_k \in\mathcal T_{n}^g$ such that
for each $(g,n)$-link  $L\subset M$ we have  $L\cong
\hat{\sigma}^{\psi_i}$ with  $\sigma\in\ker(\Omega_{g,n})$ and
$i\in\{1,\ldots,k\}$.
\end{prop}
\begin{proof}
   By result of \cite{B1}, it is possible to choose elements
    $\psi_1,\cdots,\psi_k\in\MCG(\T_g)$ such that,  each
    $\psi\in\MCG(\T_g)$ that induces  an Heegaard
    decomposition $H_g\cup_{\psi\tau}\bar H_g$   of $M$, belongs
    to the same  double coset class of $\psi_i$  in $\MCG(\T_g)$
     modulo $\X_0^g$, for a certain $i\in\{1,\ldots,k\}$.
     Now let $\psi_1,\cdots, \psi_k \in\mathcal T_{n}^g$ such that
     $\Omega_{g,n}(\psi_i)=\psi_i$, for $i=1,\ldots, k$.   Since $L$ is a
$(g,n)$-link in $M$, there exists  $\psi\in\MCG_{2n}(\T_g)$ such that
$L=\theta_{g,n}(\psi)$ and
   $ \H_g\cup_{\Omega_{g,n}(\psi)} \bar \H_g$ is a genus $g$ Heegaard
splitting for $M$. So there exist
   $\bar{\epsilon_1},\bar{\epsilon_2}\in\X^g_0$ and $i\in\{1,\ldots k\}$
such that
   $\Omega_{g,n}(\psi_i)=\bar{\epsilon_1}\Omega_{g,n}(\psi)\bar{\epsilon_2}$.
 Since    $\Omega_{g,n}$
   restricts to a surjective  homomorphism $\X_{2n}^g\to\X^g_0$, then
there exists $\epsilon_i\in\X_{2n}^g$
 such that $\Omega_{g,n}(\epsilon_i)=\bar{\epsilon_i}$, for $i=1,2$.
If we set $\sigma=\psi_i^{-1}\epsilon_1\psi\epsilon_2$, then
$\sigma\in\ker(\Omega_{g,n})$ and $L\cong \hat{\sigma}^{\psi_i}$.
\end{proof}

The Waldhausen conjecture, which has been proved in
\cite{J,J2,L}, tells us that every manifold admits  a finite
number of homeomorphism classes of irreducible genus $g$ Heegaard
splittings. So, for example,   Proposition \ref{standard} holds
whenever $g$  is the Heegaard genus of $M$.
In~\cite{CM3} it is  analyzed the case $g=n=1$.

\section{The Hilden map and the motion groups}
\label{third}
In this section we describe the connections between  Hilden braid groups
and the so-called \emph{motion groups}. We start by recalling few
definitions (see \cite{G}).

\medskip

\noindent  A \textit{motion} of a compact submanifold $N$ in a  manifold $M$ is a path $f_t$
in $Homeo_c(M)$ such that $f_0=\id_M$ and $f_1(N)=N$, where $\Homeo_c(M)$
denotes  the group of homeomorphisms of $M$ with compact support. A
motion is called \textit{stationary} if $f_t(N)=N$ for all $t\in[0,1]$. The \textit{motion
group} $\M(M,N)$ of $N$ in $M$ is the group of equivalence classes of
motion of $N$ in $M$  where two motions $f_t,g_t$ are equivalent if
$(g^{-1}f)_t$ is homotopic relative to endpoints to a stationary motion.

Notice that the motion group of $k$ points in $M$  is the braid group
$\B_k(M)$. Moreover, since each motion is equivalent to a motion that fixes
a point $*\in M-N$ , it is possible to define a homomorphism
\begin{equation}
\label{da}\M(M,N)\to \Aut(\pi_1(M-N,*))
\end{equation}
 sending an element represented by the motion $f_t$ into the automorphism
induced on $\pi_1(M-N,*)$ by $f_1$.

We are mainly interested in the case of links in 3-manifolds.  In \cite{G}
a finite set of generators for the motion groups $\mathcal M(\Sp^3,L_n)$
of the $n$-component trivial link in $\Sp^3$ is given, while a
presentation can be found in \cite{BC}.  Moreover in \cite{G2} a
presentation for the motion group of all  torus links in $\Sp^3$ is
obtained. On the contrary, there are not known examples of computations of
motion groups of links in 3-manifolds different from $\Sp^3$.

 In \cite{Hil}  was described how to construct examples of motions
of a link $L$ in $\Sp^3$ presented as the plat closure
of a braid $\sigma\in \B_{2n}(\Sp^2)$ using the elements of
$\Hil^0_n\cap\Hil^0_n(\sigma)$.
In the following Theorem we extend this result to links in 3-manifolds via
 Hilden braid groups of a surface.

\begin{thm}\label{hildenmap}
Let $\psi\in\mathcal T^g_n$ and let $\hat{\sigma}^{\psi}$ be a link in
$M_{\psi}=\H_g\cup_{\tau\psi}\bar{\H}_g$, where $\sigma\in
\ker(\Omega_{g,n})$. There exists  a  group homomorphism, that we call the
Hilden map,
 $\mathcal H_{\psi\sigma}:\Hil_n^g \cap \Hil_n^g(\psi\sigma)
\to\M(M_{\psi},\hat{\sigma}^{\phi})$, where
$\Hil_n^g(\psi\sigma)=(\psi\sigma)^{-1}\Hil_n^g\psi\sigma$.
\end{thm}
\begin{proof}
Let $\epsilon$ be a representative of an element in  $\Hil_n^g \cap
\Hil_n^g(\psi\sigma)$. By definition $\Hil_n^g\subset \ker \Omega_{g,n}$,
so there exists an isotopy $g:I\times\T_g\to T_g$  such that
$g(0,\cdot)=g_0=\textup{id}$ and $g(1,\cdot)=g_1=\epsilon$. Then
$\psi\sigma g(\psi\sigma)^{-1}$ is an isotopy between the identity and
$\psi\sigma\epsilon(\psi\sigma)^{-1}$. Moreover by hypothesis the isotopy
class of $\psi\sigma\epsilon(\psi\sigma)^{-1}$ belongs to $\Hil_n^g$ and
so extends to $\H_g$.     Since the rows of the commutative diagram
\begin{equation}
\begin{CD}
\MCG_n(\H_g)@>>>\X_{2n}^g\\
@VVV     @VV\Omega_{g,n}V\\
\MCG(\H_g)@>>>\X_{0}^g.
\end{CD}
\end{equation}  are isomorphisms,  there exist two  isotopies $f:I\times \H_g\to \H_g$
and $\bar f:I\times \H_g\to \H_g$  between the identity and an extension
of, respectively, $\epsilon$ and $\psi\sigma\epsilon(\psi\sigma)^{-1}$.
We claim that it is possible to choose $f$ and $\bar f$ such that they extend,
respectively, $g$ and $\psi\sigma g(\psi\sigma)^{-1}$. 
Indeed if $g=0$ we can use the Alexander trick to extend 
the isotopy from the boundary sphere to the 3-ball. 
If $g>0$ first we extend the isotopy  on a system of meridian discs for 
$\H_g$ not intersecting the system of arcs and then we reduce to the previous 
case by cutting along them.

We define $\mathcal{H}_{\psi\sigma}([\epsilon])=[F]$ where
$F:I\times M_{\psi}\to M_{\psi}$ is  defined by

\begin{equation*}F(t,x)F_t(x)=\left\{\begin{array}{l} f(t,x)\ \ \ \
\textup{if } x\in \H_g\\
             \bar{f}(t,x)\ \ \ \ \textup{if } x\in \bar{H}_g 
\end{array}\right.\end{equation*}
The commutativity of the following diagram ensures that $F_t$ is a
well-defined homeomorphism of $M_{\psi}$
\begin{equation}
\begin{CD}
\partial\H_g@>\tau\psi\sigma>>\partial\bar{\H}_g\\
@Vg_tVV     @VV\tau\psi\sigma g_t(\tau\psi\sigma)^{-1}V\\
\partial\H_g@>\tau\psi\sigma>>\partial\bar{\H}_g.
\end{CD}
\end{equation}
It is immediate to check  that $F_0=\textup{id}$ and
$F_1(\hat{\sigma}^{\phi})=\hat{\sigma}^{\phi}$ and so $F$ is a motion of
$\hat{\sigma}^{\phi}$ in $M_{\psi}$. Moreover the  definition  of
$\mathcal H_{\psi\sigma}$ does not depend on the homeomorphism choosen as
a representative of the element in $\Hil_n^g \cap \Hil_n^g(\psi\sigma)$:
indeed, if $\epsilon'$ is another representative,  there exists an isotopy
 between  $\epsilon$ and $\epsilon'$  fixing  $A_1\cup\cdots\cup A_n$ and
so the corresponding motions are equivalent. 

In order to prove both that the definition does not depend on the choice of the isotopies and  that 
$\mathcal{H}_{\psi\sigma}$ is a group homomorphisms we distinguish three cases. For the case of $g=0$ we refer to \cite{BC,BH,Hil}.  If $g>1$,  the statement follows from the fact that $\pi_1(\mathcal H(\T_g),\textup{id})=1$, where $\mathcal H(\T_g)$ is the group of orientation preserving homeomorphisms of $\T_g$ (see \cite{Ham2}). If $g=1$, then $\pi_1(\mathcal H(\T_g),\textup{id})=\mathbb{Z}$, see \cite{Ham}. Nevertheless, since   $\MCG(\T_1)\cong \MCG_1(\T_1)$, we can suppose that all the isotopies  that we take into consideration  fix a point; so the statement follows from  $\pi_1(\mathcal H(\T_g,P),\textup{id})=1$ (see \cite{Ham})
 
\end{proof}

In order to use the Hilden map to get informations on motion groups, it
is natural to ask if $\mathcal H_{\psi\sigma}$ is surjective and/or
injective. Clearly the answer will depend on $\psi\sigma$, that is on
the ambient manifold and on the considered link. Before  giving a
(partial) answer in the case of $\Sp^3$, let us recall the main  result
of~\cite{G}

\begin{thm}[\cite{G}]\label{gold} The homomorphism
$\mathcal(\Sp^3,L_n)\to\Aut(\pi_1(\Sp^3-L_n,*))\cong\mathbb{F}_n$ is
injective and $\mathcal(\Sp^3,L_n)$ is generated  by:
\begin{itemize}
\item[$R_i$:] turn the $i$-th circle over, corresponding to the
automorphism of $\F_n$
$$\rho_i:\ \left\{\begin{array}{l} x_i\to x_i^{-1}\\x_k\to x_k\textup{ if
}k\ne i\end{array}\right.$$
\item[$T_j$:] interchange the $j$-th and $(j+1)$-th circles,
corresponding to the automorphism of $\F_n$
$$\tau_j:\ \left\{\begin{array}{l} x_j\to x_{j+1}\\ x_{j+1}\to x_j\\x_h\to
x_h\textup{ if }h\ne j,j+1\end{array}\right.$$
\item[$A_{ik}$:] pull the $i$-th circle throught the $k$-th circle,
corresponding to the automorphism of $\F_n$
 $$\alpha_{ik}:\ \left\{\begin{array}{l} x_i\to x_kx_ix_k^{-1}\\x_h\to
x_h\textup{ if }h\ne i\end{array}\right.$$
\end{itemize}
where $j=1,\ldots,n-1$,  $i,k=1,\ldots n$ and $i\ne k$.
\end{thm}

\begin{figure}[ht]
\label{parallel}
\begin{center}
\includegraphics*[totalheight=4cm]{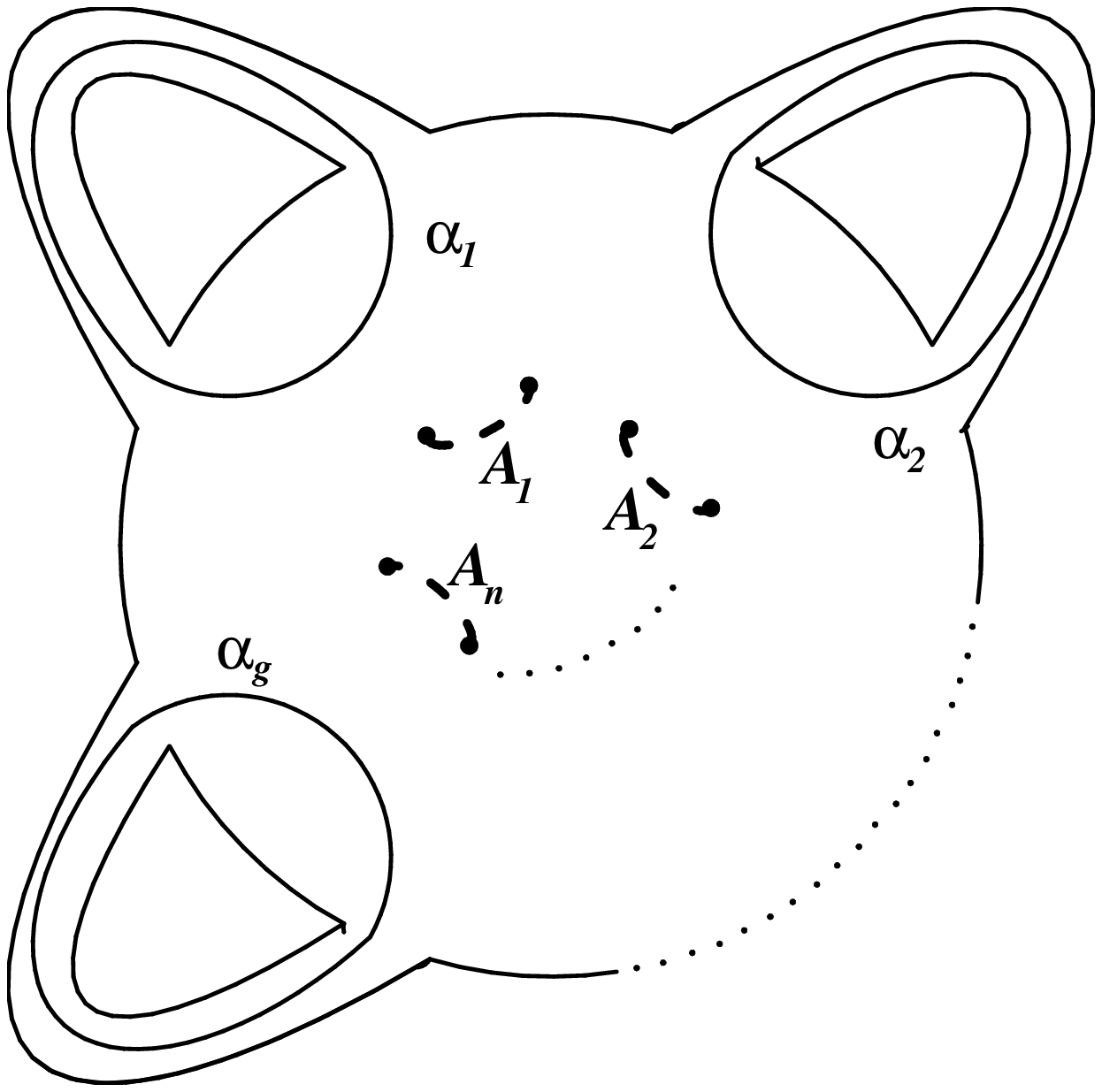}
\end{center}
\caption{An example of Hilden map.}
\end{figure}

\begin{cor}
Let $\psi\in\mathcal{T}_n^g$ be an element such that $M_{\psi}=\Sp^3$. For
example, choose $\sigma=\textup{id}$ if $g=0$ and
$\sigma=T_{\alpha_1}\cdots T_{\alpha_g}$, where $\alpha_1,\ldots\alpha_g$
denote the curves depicted in Figure~9 
 if $g>1$. The
homomorphism $\mathcal H_{\psi}: \Hil_n^g \cap \Hil_n^g(\psi)
\to\M(\Sp^3,L_n)$ is surjective. Moreover, it is injective if and only if
$(g,n)=(0,1)$.
 \end{cor}
\begin{proof}
First of all notice that each element of $\mathcal{T}_n^g$ commutes with
the following elements of $\Hil_n^g$: $\iota_i$, $\lambda_j$, $s_{i,k}$,
for $j=1,\ldots,n-1$,  $i,k=1,\ldots n$ and $i\ne k$. So all these
elements belong to $\Hil_n^g \cap \Hil_n^g(\psi)$. Moreover, $\mathcal
H_{\psi}(\iota_i)=R_i$, $\mathcal H_{\psi}(\lambda_j)=T_j$ and $\mathcal
H_{\psi}(s_{i,k})=A_{ik}$
so the surjectivity of $\mathcal H_{\psi}$ follows
by Theorem~\ref{gold}. The same holds for $\mathcal M(\Sp^3,L_1)$, since, by
Theorem~\ref{gold}, it is isomorphic to the subgroup of
$\Aut(\mathbb{F}_2)$ generated by $\rho_1$ which has clearly order two. 
On the contrary, if $(g,n)\ne(0,1)$, then $\iota_1$ has infinite order in
$\MCG_{2n}(\T_g)$  while $\rho_1$, and so $R_1$, is  an element of order
two.
\end{proof}

Using results from \cite{CM3} and \cite{G2} it would be possible to
analyze  the case of the torus links  in $\Sp^3$.
Moreover, the Hilden map could be used in order to get informations on
motion groups of links that belongs in 3-manifolds different from $\Sp^3$.

\vspace{10pt}

\noindent PAOLO BELLINGERI, Univ. Caen, CNRS UMR 6139, LMNO, Caen, 14000 (France).
Email: paolo.bellingeri@math.unicaen.fr

\vspace{5pt}
\noindent ALESSIA CATTABRIGA, Department of Mathematics - University of Bologna
Piazza di Porta S. Donato 5
40126 Bologna (Italy). Email: cattabri@dm.unibo.it

\end{document}